\documentclass{tran-l}
\usepackage{amssymb}

\usepackage{amsmath}

\input{tcilatex}
\theoremstyle{definition}
\theoremstyle{remark}
\numberwithin{equation}{section}

\input tcilatex

\begin{document}
\title[Expansion of real valued meromorphic functions...]{Expansion of real valued meromorphic functions into \textit{Fourier}
trigonometric series }
\author{\textbf{branko saric}}
\address{\textit{The Institute ''Kirilo Savic``, 11000 Belgrade, V. Stepe 51., Serbia}}
\email{bsaric@ptt.yu}
\date{May 25, 2000}
\maketitle

\begin{abstract}
In the main part of the paper, on the basis of contour integration of
complex meromorphic functions whose singularities lie onto an integration
contour, in the first step, a concept of improper integrals absolute
existence of meromorphic functions, as more general one with respect to the
concept of improper integrals convergence (existence), is introduced into
analysis. In the second step, in the case when a modulus of complex
parameter tends to infinity, an interval of improper integrals convergence
of parametric meromorphic functions is defined. In accordance with this, it
is shown that the class of real valued meromorphic functions, whose finitely
many isolated singularities lie onto a real axis segment $\left[
t_{0},\,t_{1}\right] $, may be expanded into \textit{Fourier} trigonometric
series, separately. At all points of the segment, at which the meromorphic
functions are continuous ones, the \textit{Fourier} trigonometric series is
summable and its sum is equal to the function values at those points.
Finally, that all is illustrated by two representative examples.
\end{abstract}

\section{Introduction}

From the author's viewpoint, proper attention should be paid to a class of
mathematical expressions reducing, in the boundary case, to the difference
of infinities $\infty -\infty $. Namely, in general case, the difference of
infinities $\infty -\infty $ is an indefinite expression taking any value
from the extended numeric straight-line (real axis) $R^{*}$, that is, from
the segment $\left[ -\infty ,\,+\infty \right] $ . Causality related to the
value itself of an indefinite expression is behavior of the mathematical
expression reducing to it in the boundary case. Therefore, a case of
so-called alternative numerical series is indicative one, \cite{Mi}. In
fact, since each alternative series can be expressed, according to its
definition - \textit{Definition 1, Section 2.6, Chapter 2, p. 28,} \cite{Mi}
- by difference of two positive numerical series, then in the case when each
of them definitely diverges the alternative series reduces to an indefinite
expression $\infty -\infty $. The typical representative of such a class of
series, is the alternative numerical series $\underset{k=1}{\overset{+\infty 
}{\sum }}\left( -1\right) ^{k}$. By means of \textit{Caushy's} definition of
both the sequence convergence and the numerical series convergence \cite{Mi}%
, it can be immediately shown that the alternative numerical series $%
\underset{k=1}{\overset{+\infty }{\sum }}\left( -1\right) ^{k}$,
indefinitely diverges. In other words, since in this emphasized case, the
limiting value of partial sums $\underset{k=1}{\overset{n}{\sum }}\left(
-1\right) ^{k}$ of the numerical series $\underset{k=1}{\overset{+\infty }{%
\sum }}\left( -1\right) ^{k}$ does not exist when $n\rightarrow +\infty $,
then the sum of the observed numerical series, in the \textit{Caushy's}
sense, does not exist too. Accordingly, it is natural to ask the following
questions\textit{:} How much is the sum value of this numerical series, more
exactly, is this numerical series summable? Closely related to these
questions is another\textit{:} How much is, in this emphasized instance, the
numerical value of an indefinite expression of difference of infinities $%
\infty -\infty $ ? Clearly, conceptually distinction should be made between
summation of series in the \textit{Caushy's} sense and its summability. In
the modern mathematical analysis, more exactly, in the series theory, and
from the point of view of the general convergence (summability) of numerical
series, an answer to the former questions was given by \textit{Frobenius}, 
\textit{Holder} and \textit{Cesaro} \cite{Mi}. In this paper it is presented
slightly different and more indirect answer more related to the problem of
exact determining rather than to the problem of redefining the sum itself of
indefinitely divergent series. It is more indirect because the fundamental
conclusions will be based on the results of the complex analysis theory,
more exactly, on contour integration of the functions of a complex variable.
In this case also, a proper attention will be paid to the mathematical
expressions reducing, in the boundary case, to an indefinite expression of
difference of infinities $\infty -\infty $.

\section{The main results}

\subsection{Contour integration and improper integral}

The concept of an improper integral absolute existence of a meromorphic
function $f\left( z\right) $\textit{:} $C^{1}\rightarrow C^{1}$\textit{;} $%
C\,$- is the set of complex numbers, clearly in the case when the
singularities of the function lie onto the integration contour, is based on
concept of total value ($v.t.$) of an improper integral of a meromorphic
function which is defined to be the sum of \textit{Caushy's} principal value 
\textit{(}$v.p.$\textit{)} and \textit{Jordan's} singular value \textit{(}$%
v.s.$\textit{)} of an improper integral of a meromorphic function $f\left(
z\right) $. \textit{Jordan's} singular value \textit{(}$v.s.$\textit{)} of
an improper integral of a meromorphic function is defined to be a limiting
value, as $\varepsilon \rightarrow 0^{+}$, of an integral of the function $%
f\left( z\right) $ over a certain part $\overset{\curvearrowleft }{PQ}$ of a
circular path of integration $G_{\varepsilon }$\textit{:} $G_{\varepsilon
}=\left\{ z\mid z\left( \theta \right) =\varepsilon e^{i\theta }\text{%
\textit{;} }\theta \in \left[ 0,2\pi \right] \right\} $, bypassing a
singularity of the function, where the points of the complex plane\textit{:} 
$P$ and $Q$, are intersection points of the circular contour $G_{\varepsilon
}$ and an integration contour $G$, $i\,$- denotes an imaginary unit and $e\,$%
- is a base of natural logarithm. In fact, in other words a concept of
improper integrals absolute existence of meromorphic functions generalizes
the fundamental concept of improper integrals convergence (existence).

By results of both so-called \textit{Jordan's} lemma - \textit{Theorem 1,
Subsection 3.1.4, Section 3.1, Chapter 3, p. 52}, \cite{M-K (78)} - and the
fundamental \textit{Caushy's} theorem on residues -\textit{\ Theorem 1,
Subsection 3.6.2, Section 3.6, Chapter 3, p. 226}, \cite{M-K (81)} - the sum
of \textit{Caushy's} principal value ($v.p.$) and \textit{Jordan's} singular
value ($v.s.$) of an improper integral of a meromorphic function $f\left(
z\right) $ whose only singularity (simple pole) lies onto closed integration
contour $G$, can be proved to be\footnote{{\footnotesize Symbol }$\underset{G%
}{\overset{\circlearrowleft }{\int }}${\footnotesize \ denotes an
integration over the closed contour of integration G, in this case in the
positive mathematical direction}} 
\begin{equation}
v.t.\overset{\circlearrowleft }{\underset{G}{\int }}f\left( z\right) dz=v.p.%
\overset{\circlearrowleft }{\underset{G}{\int }}f\left( z\right) dz+v.s.%
\overset{\circlearrowleft }{\underset{G_{\varepsilon }}{\int }}f\left(
z\right) dz=  \label{1}
\end{equation}
\begin{equation*}
=v.p.\overset{\circlearrowleft }{\underset{G}{\int }}f\left( z\right)
dz+\left\{ 
\begin{array}{l}
-i\alpha A \\ 
i\left( 2\pi -\alpha \right) A
\end{array}
\right. =\left\{ 
\begin{array}{l}
0 \\ 
i2\pi A
\end{array}
\right. ,
\end{equation*}
where $\alpha $ is an absolute value of an angular difference of arguments
of intersection points\textit{:} $Q$ and $P$, with respect to the point $%
z_{0}$, respectively, in the limit as $\varepsilon \rightarrow 0^{+}$, and $%
A $ is a residue of the function $f\left( z\right) $ at the point $z_{0}$,
that is\textit{:} $A=\underset{z\rightarrow z_{0}}{\lim }\left(
z-z_{0}\right) f\left( z\right) $, on condition that such limiting value
exists, \cite{M-K (78)} and \cite{St}. In this emphasized case, as
distinguished from \textit{Caushy's} principal value ($v.p.$), \textit{%
Jordan's} singular value ($v.s.$), just as well as the total value ($v.t.$),
of an improper integral of a meromorphic function $f\left( z\right) $, are
not uniquely defined ones, already they depend upon the choice of the part
of a circular path bypassing the singularity of the function.

In the case when the singularity of the meromorphic function $f\left(
z\right) $ at the point $z_{0}$ is a pole of a higher order, the conclusion
essentially differs from the preceding one. Namely, in that case, both the 
\textit{Caushy's} principal value ($v.p.$) and \textit{Jordan's} singular
value ($v.s.$) of an improper integral do not exist in the limit as $%
\varepsilon \rightarrow 0^{+}$. The improper integral reduces to an
indefinite expression of the difference of infinities $\infty -\infty $.
Thus, as for the meromorphic function $z\mapsto \frac{1}{\left(
z-z_{0}\right) ^{k}}$ , where $k\geq 2$ ($k\in N$) and $N$ is a set of
natural numbers, the improper integral along an any closed integration path
passing through the point $z_{0}$ absolutely exists and its unique total
value ($v.t.$) is identically zero. This is in agreement with both the
general \textit{Caushy-Goursat's} integral theorem - \textit{Theorem 1,
Subsection 3.5.2, Section 3.5, Chapter 3, p. 203}, \cite{M-K (81)} - and the
results of contour integration of rational functions - \textit{Subsubsection
3.1.2.3, Subsection 3.1.2. Section 3.1, Chapter 3, p. 46}, \cite{M-K (78)}.
Now then, in this emphasized case, a sum of values of integrals of
meromorphic function $f\left( z\right) $ over a part of any integration path 
$G$ between intersection points as well as over the part $\overset{%
\curvearrowleft }{PQ}$ of circular path $G_{\varepsilon }$, is identically
zero for each $\varepsilon $. Since choice of an arc path\textit{:} $%
\overset{\curvearrowleft }{PQ}$ or $\overset{\curvearrowright }{PQ}$,
bypassing a singularity of the function is arbitrary, the above-mentioned
unique sum remains zero in the limit as $\varepsilon \rightarrow 0^{+}$.

\subsubsection{Example}

Let the part of an integration path between points\textit{:} $P$ and $Q$, be
a part of circumference of a circle centred at the orgin and of radius $a$%
\textit{:} $a\in R_{+}^{1}$ ($R_{+}^{1}$- is a set of positive real numbers)
and $a$ be also a singularity of the function $f\left( z\right) $\textit{:} $%
f\left( z\right) =\frac{1}{\left( z-a\right) ^{k^{*}}}$, $k^{*}\geq 2$ and $%
k^{*}\in N$. Then, an integral of the function $f\left( z\right) $, over the
part of a circular integration path $G$ from the point $P$ to the point $Q$
which dose not contain the singularity $a$, reduces to the integral 
\begin{equation}
a\overset{2\pi -\alpha ^{*}}{\underset{\alpha ^{*}}{\int }}\frac{e^{i\theta
^{*}}}{\left( ae^{i\theta ^{*}}-a\right) ^{k^{*}}}d\theta ^{*}=  \label{2}
\end{equation}
\begin{equation*}
=-\frac{\left( -i\right) ^{k}}{\left( 2a\right) ^{k+1}}\left[ \overset{\pi
-\alpha }{\underset{\alpha }{\int }}\frac{\cos \left( k\theta \right) }{%
\left( \sin \theta \right) ^{k+2}}d\theta -i\overset{\pi -\alpha }{\underset{%
\alpha }{\int }}\frac{\sin \left( k\theta \right) }{\left( \sin \theta
\right) ^{k+2}}d\theta \right] \text{,}
\end{equation*}
where\thinspace $\alpha ^{*}$ is an argument of the point $P$ with respect
to the origin ($\alpha =\frac{\alpha ^{*}}{2}$), $\theta =\frac{\theta ^{*}}{%
2}$ and $k=k^{*}-2$ ($k\in N_{0}$\textit{;} $N_{0}=\left\{ 0,1,2,...\right\} 
$).

As results of partial integration the following integral dependencies are
obtained\textit{:} 
\begin{equation}
i\left( 1+k\right) \overset{\pi -\alpha }{\underset{\alpha }{\int }}\frac{%
2ae^{2i\theta }}{\left( ae^{2i\theta }-a\right) ^{k+2}}d\theta =\left( \frac{%
-i}{2a}\right) ^{k+1}\left[ -\frac{\cos \left( k\theta \right) }{\left( \sin
\theta \right) ^{k}}\cot \theta \left| _{\alpha }^{\pi -\alpha }\right.
-\right.  \label{3}
\end{equation}
\begin{equation*}
\left. -k\overset{\pi -\alpha }{\underset{\alpha }{\int }}\frac{\sin \left[
\left( k-1\right) \theta \right] }{\left( \sin \theta \right) ^{k+1}}d\theta
-i\left( k+1\right) \overset{\pi -\alpha }{\underset{\alpha }{\int }}\frac{%
\sin \left( k\theta \right) }{\left( \sin \theta \right) ^{k+2}}d\theta
\right]
\end{equation*}
and 
\begin{equation}
\left( k+2\right) \overset{\pi -\alpha }{\underset{\alpha }{\int }}\frac{%
\sin \left[ \left( k+1\right) \theta \right] }{\left( \sin \theta \right)
^{k+3}}d\theta =-2\frac{\cos \left( k\theta \right) }{\left( \sin \theta
\right) ^{k}}\cot \theta \left| _{\alpha }^{\pi -\alpha }\right. -  \label{4}
\end{equation}
\begin{equation*}
-2k\overset{\pi -\alpha }{\underset{\alpha }{\int }}\frac{\sin \left[ \left(
k-1\right) \theta \right] }{\left( \sin \theta \right) ^{k+1}}d\theta -\frac{%
\sin \left( k\theta \right) }{\left( \sin \theta \right) ^{k+1}}\frac{1}{%
\sin \theta }\left| _{\alpha }^{\pi -\alpha }\right. \text{.}
\end{equation*}

As, in this emphasized case, $2a\sin \alpha =\varepsilon $ a by-pass
integral value is\textit{:} 
\begin{equation}
\varepsilon ^{-\left( k+1\right) }\underset{-\left( \frac{\pi }{2}+\alpha
\right) }{\overset{\frac{\pi }{2}+\alpha }{\int }}e^{-i\left( k+1\right)
\theta }d\theta =\left\{ 
\begin{array}{l}
\frac{\left( -1\right) ^{n}}{n}\frac{\sin \left( 2n\alpha \right) }{\left(
2a\sin \alpha \right) ^{2n}}\text{\textit{;} }k=2n-1 \\ 
\frac{2\left( -1\right) ^{n}}{2n+1}\frac{\cos \left[ \left( 2n+1\right)
\alpha \right] }{\left( 2a\sin \alpha \right) ^{2n+1}}\text{\textit{;} }k=2n
\end{array}
\right. \text{, }n\in N\text{.}  \label{5}
\end{equation}

By the comparative analysis of preceding equalities it can be easily shown
that if the sum of the integral value of the function $f\left( z\right) $
over the circular integration contour $G$ from the point $P$ to the point $Q$
and the by-pass integral value, for arbitrarily chosen $\alpha $ and $k=2n-1$
respectively, is identically zero, just as well as the value of the integral 
$\overset{\pi -\alpha }{\underset{\alpha }{\int }}\frac{\sin \left[ \left(
k-1\right) \theta \right] }{\left( \sin \theta \right) ^{k+1}}d\theta $ for $%
k=2n-1$, see the equation (\ref{4}), then for $k=2n-1$ and $n\in N$ it holds 
\begin{equation}
\overset{\pi -\alpha }{\underset{\alpha }{\int }}\frac{\sin \left( k\theta
\right) }{\left( \sin \theta \right) ^{k+2}}d\theta =\frac{1}{n}\frac{\sin
\left( 2n\alpha \right) }{\left( \sin \alpha \right) ^{2n}}\text{.}
\label{6}
\end{equation}

Since the second integral on the right-hand side of the equation (\ref{3})
is identically zero for $k=2n$, then in view of the preceding result (\ref{6}%
) it follows that 
\begin{equation}
\overset{\pi -\alpha }{\underset{\alpha }{\int }}\frac{2ae^{2i\theta }}{%
\left( ae^{2i\theta }-a\right) ^{2\left( n+1\right) }}d\theta =\frac{2\left(
-1\right) ^{n}}{2n+1}\left\{ \frac{-\cos \left[ \left( 2n+1\right) \alpha
\right] }{\left( 2a\sin \alpha \right) ^{2n+1}}\right\} \text{, }n\in N\text{%
.}  \label{7}
\end{equation}

If one considers the fact that, for a corresponding natural number $n$ ($%
n\in N$), the sum of functional expressions on right-hand sides of relations%
\textit{:} (\ref{5}) (for $k=2n$) and (\ref{7}), is identically zero for
arbitrarily chosen $\alpha $, then the sum of integrals on the left-hand
sides of these equations is identically zero too. Thus, the improper
integral of the function $f\left( z\right) $ along the circular contour of
integration $G$ absolutely exists and its total value, as limiting value of
a sum of integrals (\ref{2}) and (\ref{5}) as $\alpha \rightarrow 0^{+}$, is
identically zero for each $k^{*}\geq 2$ and $k^{*}\in N$, and what has been
just proved by method of a mathematical induction\footnote{{\footnotesize %
For }$k=0${\footnotesize \ and }$k=1${\footnotesize \ i.e. }$k^{*}=2$%
{\footnotesize \ and }$k^{*}=3${\footnotesize , on the basis of the relation
(\ref{3}) it holds }$\overset{\pi -\alpha }{\underset{\alpha }{\int }}\frac{%
2ae^{2i\theta }}{\left( ae^{2i\theta }-a\right) ^{k+2}}d\theta =-\frac{\cot
\alpha }{a}${\footnotesize \ and }$\overset{\pi -\alpha }{\underset{\alpha }{%
\int }}\frac{2ae^{2i\theta }}{\left( ae^{2i\theta }-a\right) ^{k+2}}d\theta =%
\frac{\cot \alpha }{2a^{2}}${\footnotesize , respectively. From the relation
(\ref{5}), the values of by-pass integrals, in these emphasized cases, are
equal to\textit{:} }$\frac{1}{\varepsilon }\underset{-\left( \frac{\pi }{2}%
+\alpha \right) }{\overset{\frac{\pi }{2}+\alpha }{\int }}e^{-i\theta
}d\theta =\frac{\cot \alpha }{a}${\footnotesize \ and }$\frac{1}{\varepsilon
^{2}}\underset{-\left( \frac{\pi }{2}+\alpha \right) }{\overset{\frac{\pi }{2%
}+\alpha }{\int }}e^{-2i\theta }d\theta =-\frac{\cot \alpha }{2a^{2}}$%
{\footnotesize , respectively. Accordingly, the total value of improper
integral }$i\underset{o}{\overset{2\pi }{\int }}\frac{ae^{i\theta ^{*}}}{%
\left( ae^{i\theta ^{*}}-a\right) ^{k^{*}}}d\theta ^{*}${\footnotesize , is
identically zero, for both }$k^{*}=2${\footnotesize \ and }$k^{*}=3$%
{\footnotesize .}}.$\blacktriangledown $

\begin{description}
\item[Comment]  As it has been just illustrated by the previous example, the
concept of an improper integral absolute existence is more general concept
with respect to the concept of an improper integral convergence. This is in
connection with indefinite expression of the difference of infinities $%
\infty -\infty $ to which the improper integral value is reduced in a
boundary case. Namely, independently from the fact that the improper
integral absolutely exists, in some of the concrete cases its principal
value does not exist. Hence, by introducing a by-pass integral into the
analysis the concept itself of improper integral convergence (existence) is
generalized to the concept of improper integral absolute existence.$%
\blacktriangledown $
\end{description}

\subsection{Fourier trigonometric series of real valued meromorphic functions
}

\subsubsection{An analysis of an idea}

Without loss of the generality, one may assume that a complex meromorphic
function $g\left( z,t\right) $, where the variable $t$ is independent one
with respect to the complex variable $z$, has infinitely but a count of many
simple poles\textit{:} $a_{1},a_{2},...$ onto the imaginary axis of the
complex plane $C^{1}$. In that emphasized case, there exists a sequence of
circular contours of integration $G_{r}$, centred at the origin and of
radius $r$, such that onto theirs boundaries there are no singularities of
the function $g\left( z,t\right) $. Hence, by the fundamental \textit{%
Cauchy's} theorem on residues the sequence of the partial sums can be formed 
\begin{equation}
\overset{n}{\underset{k=1}{\sum }}A_{k}\left( t\right) =\frac{1}{2\pi i}%
\overset{\circlearrowleft }{\underset{G_{r}}{\int }}g\left( z,t\right) dz%
\text{,}  \label{8}
\end{equation}
where $A_{k}\left( t\right) $\textit{:} $A_{k}\left( t\right) =\underset{%
z=a_{k}}{Res}g\left( z,t\right) $, are residues of the function $g\left(
z,t\right) $ at the points\textit{:} $z=a_{k}$, $k=1,2,...,n$.

On the one hand, on the basis of the second \textit{Jordan's} lemma - 
\textit{Theorem 2, Subsection 3.1.4, Section 3.1, Chapter 3, p. 52}, \cite
{M-K (78)} - if there exists an unique limiting value\textit{:} $\underset{%
\left| z\right| \rightarrow +\infty }{\lim }\left[ zg\left( z,t\right)
\right] $\textit{;} for each $z\in C^{1}$, then the sequence of the partial
sums $\overset{n}{\underset{k=1}{\sum }}A_{k}\left( t\right) $ converges, in
other words there exists a sum of the infinite functional series $\overset{%
+\infty }{\underset{k=1}{\sum }}A_{k}\left( t\right) $ in the \textit{%
Cauchy's} sense\textit{:} 
\begin{equation}
\overset{+\infty }{\underset{k=1}{\sum }}A_{k}\left( t\right) =-\underset{%
\left| z\right| =+\infty }{Res}g\left( z,t\right) \text{,}  \label{9}
\end{equation}
where $2\pi i\underset{\left| z\right| =+\infty }{Res}g\left( z,t\right) =%
\underset{r\rightarrow +\infty }{\lim }\overset{\circlearrowright }{%
\underset{G_{r}}{\int }}g\left( z,t\right) dz=-2\pi i\underset{\left|
z\right| \rightarrow +\infty }{\lim }\left[ zg\left( z,t\right) \right] $.

However, on the other hand, by the same \textit{Jordan's} lemma, if there
exists no an above mentioned limiting value\textit{:} $\underset{\left|
z\right| \rightarrow +\infty }{\lim }\left[ zg\left( z,t\right) \right] $%
\textit{;} for each $z\in C^{1}$, already there exist only partial limiting
values\textit{:} $\underset{\left| z\right| \rightarrow +\infty }{\lim }%
\left[ zg\left( z,t\right) \right] $\textit{; }$\func{Re}z>0$ and $\underset{%
\left| z\right| \rightarrow +\infty }{\lim }\left[ zg\left( z,t\right)
\right] $\textit{;} $\func{Re}z<0$, then there exists a infinite sum of the
residues of the function $g\left( z,t\right) $ that is equal to the limiting
sum of integral values\textit{:} 
\begin{equation}
\overset{+\infty }{\underset{k=1}{\sum }}A_{k}\left( t\right) =\frac{1}{2\pi
i}\underset{r\rightarrow +\infty }{\lim }\left[ \overset{\curvearrowleft }{%
\underset{G_{r}^{R}}{\int }}g\left( z,t\right) dz+\overset{\curvearrowleft }{%
\underset{G_{r}^{L}}{\int }}g\left( z,t\right) dz\right] \text{,}  \label{10}
\end{equation}
where the integral paths\textit{:} $G_{r}^{R}=\left\{ z\mid z\left( \theta
\right) =re^{i\theta }\text{\textit{;} }\theta \in \left[ -\frac{\pi }{2}%
+\delta \left( r\right) ,\frac{\pi }{2}-\delta \left( r\right) \right]
\right\} $ and $G_{r}^{L}=\left\{ z\mid z\left( \theta \right) =re^{i\theta }%
\text{\textit{;} }\theta \in \left[ \frac{\pi }{2}+\delta \left( r\right) ,%
\frac{3\pi }{2}-\delta \left( r\right) \right] \right\} $, are arc parts of
the circular path of integration $G_{r}$ in the right-hand and left-hand
half-plane of the complex plane $C^{1}$, respectively, and an arbitrary
angular function $\delta \left( r\right) $, which is of sufficiently small
real positive values for any positive values of $r$, satisfies the condition%
\textit{:} $\underset{r\rightarrow +\infty }{\lim }\delta \left( r\right) =0$%
. In other words, although in this emphasized case there exists no a sum of
the infinite functional series $\overset{+\infty }{\underset{k=1}{\sum }}%
A_{k}\left( t\right) $ in the \textit{Cauchy's} sense, this infinite
functional series is summable\textit{.} Note that in this case too\textit{:} 
$\overset{+\infty }{\underset{k=1}{\sum }}A_{k}\left( t\right) =-\underset{%
\left| z\right| =+\infty }{Res}g\left( z,t\right) $, where $2\pi i\underset{%
\left| z\right| =+\infty }{Res}g\left( z,t\right) =-\underset{r\rightarrow
+\infty }{\lim }\left[ \overset{\curvearrowleft }{\underset{G_{r}^{R}}{\int }%
}g\left( z,t\right) dz+\overset{\curvearrowleft }{\underset{G_{r}^{L}}{\int }%
}g\left( z,t\right) dz\right] $.

\subsubsection{\textit{Cauchy's} formula}

As it is well-known, during the deriving \textit{Caushy's} formula for
expansion of real valued functions into an infinite functional series, \cite
{Ca} and \cite{Pi} (taken over from \cite{M-K (78)}) - \textit{Formula (9),
Subsection 4.6.2, Section 4.6, Chapter 4, p. 94}, \cite{M-K (78)} - in a
real axis interval $\left( t_{0},t_{1}\right) $ 
\begin{equation}
f\left( t\right) =\overset{+\infty }{\underset{k=1}{\sum }}\frac{w\left(
a_{k}\right) }{\frac{dq\left( z\right) }{dz}\left| _{z=a_{k}}\right. }%
\underset{t_{0}}{\overset{t_{1}}{\int }}e^{a_{k}\left( t-\tau \right)
}f\left( \tau \right) d\tau \text{\textit{;} }t\neq t_{si}\text{,}
\label{11}
\end{equation}
where $t_{si}$ are break points of the function $t\mapsto f\left( t\right) $
in $\left( t_{0},t_{1}\right) $, the conditions for existence of finite
limiting values of the functional expressions\textit{:} 
\begin{equation}
\underset{\left| z\right| \rightarrow +\infty }{\lim }\left[ z\frac{p\left(
z\right) }{q\left( z\right) }\underset{t_{0}}{\overset{t}{\int }}e^{z\left(
t-\tau \right) }f\left( \tau \right) d\tau -z\frac{p\left( -z\right) }{%
q\left( -z\right) }\underset{t_{0}}{\overset{t}{\int }}e^{-z\left( t-\tau
\right) }f\left( \tau \right) d\tau \right] \text{\textit{;}}  \label{12}
\end{equation}
\begin{equation}
\underset{\left| z\right| \rightarrow +\infty }{\lim }\left[ z\frac{w\left(
z\right) }{q\left( z\right) }\underset{t}{\overset{t_{1}}{\int }}e^{z\left(
t-\tau \right) }f\left( \tau \right) d\tau -z\frac{w\left( -z\right) }{%
q\left( -z\right) }\underset{t}{\overset{t_{1}}{\int }}e^{-z\left( t-\tau
\right) }f\left( \tau \right) d\tau \right] \text{,}  \label{13}
\end{equation}
which have to be satisfied by the function $t\mapsto f\left( t\right) $ in $%
\left( t_{0},t_{1}\right) $, are of the most importance.

Namely, let $g_{1}\left( z,t\right) =\frac{p\left( z\right) }{q\left(
z\right) }\underset{t_{0}}{\overset{t}{\int }}e^{z\left( t-\tau \right)
}f\left( \tau \right) d\tau $ and $g_{2}\left( z,t\right) =\frac{w\left(
z\right) }{q\left( z\right) }\underset{t}{\overset{t_{1}}{\int }}e^{z\left(
t-\tau \right) }f\left( \tau \right) d\tau $, where an analytic function $%
q\left( z\right) $\textit{:} $q\left( z\right) =p\left( z\right) +w\left(
z\right) $, has infinitely but a count of many simple poles\textit{:} $%
a_{1},a_{2},...$ onto the imaginary axis as singularities. If under an
assumption that\textit{:} $\underset{\left| z\right| \rightarrow +\infty }{%
\lim }\frac{p\left( z\right) }{q\left( z\right) }e^{z\left( t-t_{0}\right)
}=0$ and $\underset{\left| z\right| \rightarrow +\infty }{\lim }\frac{%
p\left( -z\right) }{q\left( -z\right) }=1$ as well as $\underset{\left|
z\right| \rightarrow +\infty }{\lim }\frac{w\left( -z\right) }{q\left(
-z\right) }e^{z\left( t_{1}-t\right) }=0$ and $\underset{\left| z\right|
\rightarrow +\infty }{\lim }\frac{w\left( z\right) }{q\left( z\right) }=1$,
the following functional expressions\textit{:} $\underset{\left| z\right|
\rightarrow +\infty }{\lim }z\underset{t_{0}}{\overset{t}{\int }}e^{-z\left(
\tau -t_{0}\right) }f\left( \tau \right) d\tau $ and$\underset{\left|
z\right| \rightarrow +\infty }{\lim }z\underset{t}{\overset{t_{1}}{\int }}%
e^{-z\left( t_{1}-\tau \right) }f\left( \tau \right) d\tau $, converge for
each $\func{Re}z>0$ and $t\in $ $\left( t_{0},t_{1}\right) $, just as well
as functional expressions\textit{:} $\underset{\left| z\right| \rightarrow
+\infty }{\lim }z\underset{t_{0}}{\overset{t}{\int }}e^{-z\left( t-\tau
\right) }f\left( \tau \right) d\tau $ and$\underset{\left| z\right|
\rightarrow +\infty }{\lim }z\underset{t}{\overset{t_{1}}{\int }}e^{-z\left(
\tau -t\right) }f\left( \tau \right) d\tau $, in such a way that\textit{:} $%
\underset{\left| z\right| \rightarrow +\infty }{\lim }z\underset{t_{0}}{%
\overset{t}{\int }}e^{-z\left( t-\tau \right) }f\left( \tau \right) d\tau
=f\left( t\right) $ and$\underset{\left| z\right| \rightarrow +\infty }{\lim 
}z\underset{t}{\overset{t_{1}}{\int }}e^{-z\left( \tau -t\right) }f\left(
\tau \right) d\tau =f\left( t\right) $\textit{;} $t\neq t_{si}$, then on the
basis of previous analyzed idea it can be proved that at all points of $%
\left( t_{0},t_{1}\right) $, at which a function $t\mapsto f\left( t\right) $
is continuous, there exists a sum of an infinite functional series on the
right-hand side of the equation (\ref{11}) which is just equal to the
function values at those points, more exactly, in the general case, \textit{%
Cauchy's} infinite functional series of the function $t\mapsto f\left(
t\right) $ is summable.

If a function $t\mapsto f\left( t\right) $ satisfies the general well-known 
\textit{Dirichlet's} conditions in $\left( t_{0},t_{1}\right) $, then
partial sums of \textit{Cauchy's} infinite functional series of the function 
$t\mapsto f\left( t\right) $ at all points of $\left( t_{0},t_{1}\right) $,
at which the function $t\mapsto f\left( t\right) $ is continuous, converge
to the function values at those points \cite{M-K (78)}. At the break points $%
t_{si}$ of the function $t\mapsto f\left( t\right) $ in $\left(
t_{0},t_{1}\right) $ the partial sums of \textit{Cauchy's} infinite
functional series of the function $t\mapsto f\left( t\right) $ converge to
the following functional values \cite{M-K (78)} 
\begin{equation}
\frac{1}{2}\left[ \underset{\varepsilon \rightarrow 0^{+}}{\lim }f\left(
t_{si}+\varepsilon \right) +\underset{\eta \rightarrow 0^{+}}{\lim }f\left(
t_{si}-\eta \right) \right] \text{.}  \label{14}
\end{equation}

At the extreme points of the segment $\left[ t_{0},t_{1}\right] $\textit{:} $%
t_{0}$ and $t_{1}$, at which a function $t\mapsto f\left( t\right) $ is
continuous on one's right and left respectively, the sum of \textit{Cauchy's}
infinite functional series of the function $t\mapsto f\left( t\right) $ is
equal to the following functional value 
\begin{equation}
\frac{1}{2}\left[ \underset{\varepsilon \rightarrow 0^{+}}{\lim }f\left(
t_{0}+\varepsilon \right) +\underset{\eta \rightarrow 0^{+}}{\lim }f\left(
t_{1}-\eta \right) \right] \text{.}  \label{15}
\end{equation}

\subsubsection{Interval of improper integrals convergence of real valued
functions}

Let $\nu \mapsto f\left( \nu \right) $, be an analytic function of complex
variable $\nu $ on some neighborhood $V_{0}$ of the point $\nu =0$ at which
the function $\nu \mapsto f\left( \nu \right) $ has a pole of arbitrary
order as a singularity. A function $\nu \mapsto f\left( \nu \right)
e^{-z\left( t-\nu \right) }$, where $z\in C^{1}$ is a complex parameter and $%
t\in R_{+}^{1}$ ($\func{Re}\nu =t$) is fixed point belonging to the
neighborhood $V_{0}$, is parametric analytic function on $V_{0}$. Further, a
smooth one-one mapping $\nu \left( \theta \right) $\textit{:} $%
R^{1}\rightarrow C^{1}$ ($\nu \left( \theta \right) =\frac{t-t_{0}}{2}%
e^{i\theta }+\frac{t_{0}+t}{2}$) of the real axis segment $\left[ -\pi
,0\right] $ ($\theta \in \left[ -\pi ,0\right] $) onto the set of complex
points $\nu $ of the complex plane $C^{1}$ is defined. A fixed point $\func{%
Re}\nu =t_{0}$ ($t_{0}<0$) also belongs to the neighborhood $V_{0}$ of the
zero point $\nu =0$. An arbitrary $n$-division $P_{n}$\textit{:} $%
P_{n}=\left\{ \theta _{0}=-\pi ,\theta _{1},...,\theta _{i},...,\theta
_{n}=0\right\} $, where $n\in N$, is one of all possible $n$-divisions of
the segment $\left[ -\pi ,0\right] $.

Accordingly, since a complex function $f\left[ \nu \left( \theta \right)
\right] e^{-z\left[ t-\nu \left( \theta \right) \right] }$ of a real
variable $\theta $ is a regular that in the segment $\left[ -\pi ,0\right] $%
, and in those circumstances its both a real and an imaginary part satisfies
all conditions of \textit{Langrange's} mean value theorem of the
differential calculus in the segment $\left[ -\pi ,0\right] $, then, for
each partial segment $\left[ \theta _{i-1},\theta _{i}\right] $ of the
segment $\left[ -\pi ,0\right] $, it holds 
\begin{equation}
\func{Re}\left\{ \left\{ \frac{d}{d\theta }\left\{ f\left[ \nu \left( \theta
\right) \right] e^{-z\left[ t-\nu \left( \theta \right) \right] }\right\}
\right\} \left| _{\theta =\theta _{i}^{*}}\right. \right\} =  \label{16}
\end{equation}
\begin{equation*}
=\frac{\func{Re}\left\{ f\left[ \nu \left( \theta _{i}\right) \right]
e^{-z\left[ t-\nu \left( \theta _{i}\right) \right] }\right\} -\func{Re}%
\left\{ f\left[ \nu \left( \theta _{i-1}\right) \right] e^{-z\left[ t-\nu
\left( \theta _{i-1}\right) \right] }\right\} }{\theta _{i}-\theta _{i-1}}%
\text{\textit{;}}
\end{equation*}
\begin{equation}
\func{Im}\left\{ \left\{ \frac{d}{d\theta }\left\{ f\left[ \nu \left( \theta
\right) \right] e^{-z\left[ t-\nu \left( \theta \right) \right] }\right\}
\right\} \left| _{\theta =\theta _{i}^{**}}\right. \right\} =  \label{17}
\end{equation}
\begin{equation*}
=\frac{\func{Im}\left\{ f\left[ \nu \left( \theta _{i}\right) \right]
e^{-z\left[ t-\nu \left( \theta _{i}\right) \right] }\right\} -\func{Im}%
\left\{ f\left[ \nu \left( \theta _{i-1}\right) \right] e^{-z\left[ t-\nu
\left( \theta _{i-1}\right) \right] }\right\} }{\theta _{i}-\theta _{i-1}}%
\text{,}
\end{equation*}
where $\left\{ \theta _{i}^{*},\theta _{i}^{**}\right\} \in \left[ \theta
_{i-1},\theta _{i}\right] $.

By virtue of (\ref{16}) and (\ref{17}), it is possible to form the integral
sums 
\begin{equation}
\overset{n}{\underset{i=1}{\sum }}\func{Re}\left\{ \left\{ \frac{d}{d\theta }%
\left\{ f\left[ \nu \left( \theta \right) \right] e^{-z\left[ t-\nu \left(
\theta \right) \right] }\right\} \right\} \left| _{\theta =\theta
_{i}^{*}}\right. \right\} \left( \theta _{i}-\theta _{i-1}\right) =
\label{18}
\end{equation}
\begin{equation*}
=\func{Re}\left[ f\left( t\right) \right] -\func{Re}\left[ f\left(
t_{0}\right) e^{-z\left( t-t_{0}\right) }\right] \text{\textit{;}}
\end{equation*}
\begin{equation}
\overset{n}{\underset{i=1}{\sum }}\func{Im}\left\{ \left\{ \frac{d}{d\theta }%
\left\{ f\left[ \nu \left( \theta \right) \right] e^{-z\left[ t-\nu \left(
\theta \right) \right] }\right\} \right\} \left| _{\theta =\theta
_{i}^{**}}\right. \right\} \left( \theta _{i}-\theta _{i-1}\right) =
\label{19}
\end{equation}
\begin{equation*}
=\func{Im}\left[ f\left( t\right) \right] -\func{Im}\left[ f\left(
t_{0}\right) e^{-z\left( t-t_{0}\right) }\right] \text{,}
\end{equation*}
that is, after the performed differentiation 
\begin{equation}
\overset{n}{\underset{i=1}{\sum }}\func{Re}\left\{ \left\{ \frac{d}{d\theta }%
f\left[ \nu \left( \theta \right) \right] e^{-z\left[ t-\nu \left( \theta
\right) \right] }\right\} \left| _{\theta =\theta _{i}^{*}}\right. \right\}
\left( \theta _{i}-\theta _{i-1}\right) +  \label{20}
\end{equation}
\begin{equation*}
+\overset{n}{\underset{i=1}{\sum }}\func{Re}\left\{ \left\{ zf\left[ \nu
\left( \theta \right) \right] e^{-z\left[ t-\nu \left( \theta \right)
\right] }\right\} \left| _{\theta =\theta _{i}^{*}}\right. \right\} \left(
\theta _{i}-\theta _{i-1}\right) =
\end{equation*}
\begin{equation*}
=\func{Re}\left[ f\left( t\right) \right] -\func{Re}\left[ f\left(
t_{0}\right) e^{-z\left( t-t_{0}\right) }\right] \text{\textit{;}}
\end{equation*}
\begin{equation}
\overset{n}{\underset{i=1}{\sum }}\func{Im}\left\{ \left\{ \frac{d}{d\theta }%
f\left[ \nu \left( \theta \right) \right] e^{-z\left[ t-\nu \left( \theta
\right) \right] }\right\} \left| _{\theta =\theta _{i}^{**}}\right. \right\}
\left( \theta _{i}-\theta _{i-1}\right) +  \label{21}
\end{equation}
\begin{equation*}
+\overset{n}{\underset{i=1}{\sum }}\func{Im}\left\{ \left\{ zf\left[ \nu
\left( \theta \right) \right] e^{-z\left[ t-\nu \left( \theta \right)
\right] }\right\} \left| _{\theta =\theta _{i}^{**}}\right. \right\} \left(
\theta _{i}-\theta _{i-1}\right) =
\end{equation*}
\begin{equation*}
=\func{Im}\left[ f\left( t\right) \right] -\func{Im}\left[ f\left(
t_{0}\right) e^{-z\left( t-t_{0}\right) }\right] \text{.}
\end{equation*}

Having in view the fact that $\nu \left( \theta \right) =\frac{t-t_{0}}{2}%
e^{i\theta }+\frac{t_{0}+t}{2}$, it is possible to define an interval (a
semi-interval) of a change of the argument $\varphi $: $\varphi \in R^{1}$,
of the complex parameter $z$, for which it holds: $\underset{\left| z\right|
\rightarrow +\infty }{\lim }e^{-z\left[ t-\nu \left( \theta \right) \right]
}=0$\textit{;} $\theta \in \left( -\pi ,0\right) $. Namely, since for $\func{%
Re}z\neq 0$%
\begin{equation}
z\left[ t-\nu \left( \theta \right) \right] =z\left( \frac{t-t_{0}}{2}%
\right) \left( 1-e^{i\theta }\right) =  \label{22}
\end{equation}
\begin{equation*}
=\left( \frac{t-t_{0}}{2}\right) \left[ \func{Re}z\left( 1-\cos \theta
\right) \right] \left( 1+\frac{\func{Im}z}{\func{Re}z}\frac{\sin \theta }{%
1-\cos \theta }\right) +
\end{equation*}
\begin{equation*}
+i\left[ \func{Im}z\left( 1-\cos \theta \right) -\func{Re}z\sin \theta
\right] \text{,}
\end{equation*}
then, if the condition 
\begin{equation}
\func{Re}z\left( 1-\cos \theta \right) \left( 1+\frac{\func{Im}z}{\func{Re}z}%
\frac{\sin \theta }{1-\cos \theta }\right) >0  \label{23}
\end{equation}
is satisfied, it follows that $\underset{\left| z\right| \rightarrow +\infty 
}{\lim }e^{-z\left[ t-\nu \left( \theta \right) \right] }=0$\textit{;} $%
\theta \in \left( -\pi ,0\right) $.

In view of the fact that $\frac{\sin \theta }{1-\cos \theta }=\cot \frac{%
\theta }{2}=\left( \tan \frac{\theta }{2}\right) ^{-1}$ and $\frac{\func{Im}z%
}{\func{Re}z}=\tan \varphi $, the condition (\ref{23}) is satisfied if and
only if $\varphi \in \left( -\frac{\pi }{2},0\right] $, in other words for
each $\varphi \in \left( -\frac{\pi }{2},0\right] $ it holds $\underset{%
\left| z\right| \rightarrow +\infty }{\lim }e^{-z\left[ t-\nu \left( \theta
\right) \right] }=0$\textit{;} $\theta \in \left( -\pi ,0\right) $.

Hence, and on the basis of derived relations\textit{:} (\ref{20}) and (\ref
{21}), in the limit as $n\rightarrow +\infty $, more exactly, when a maximum
partial segment $\left[ \theta _{i-1},\theta _{i}\right] $ of the segment $%
\left[ -\pi ,0\right] $ vanishes, the condition $\varphi \in \left( -\frac{%
\pi }{2},0\right] $ becomes a condition of convergence of a parametric
contour integral 
\begin{equation}
\underset{\left| z\right| \rightarrow +\infty }{\lim }z\underset{G}{\int }%
f\left( \nu \right) e^{-z\left( t-\nu \right) }d\nu =f\left( t\right) \text{,%
}  \label{24}
\end{equation}
where $G=\left\{ \nu \mid \text{ }\nu =\frac{t-t_{0}}{2}e^{i\theta }+\frac{%
t_{0}+t}{2}\right\} $ and $\theta \in \left[ -\pi ,0\right] $.

In other words, the semi-interval $\left( -\frac{\pi }{2},0\right] $ is a
semi-interval of a complex parametric contour integral convergence of the
function $f\left( \nu \right) e^{-z\left( t-\nu \right) }$ along the given
contour of integration $G$, in the limit as $\left| z\right| \rightarrow
+\infty $.

In the next step, instead of the above integration path $G$, a complex plane
curve $G^{*}$ consisting of parts of real axis defined by segments\textit{:} 
$\left[ t_{0},-\varepsilon \right] $ and $\left[ \varepsilon ,t\right] $ ($%
\varepsilon \in R_{+}^{1}$) as well as of a part of a circular path defined
by a smooth one-one mapping $\nu \left( \theta \right) $\textit{:} $%
R^{1}\rightarrow C^{1}$ ($\nu \left( \theta \right) =\varepsilon e^{i\theta
} $) of the segment $\left[ -\pi ,0\right] $ of a real axis onto a set of
points of the complex plane, is taken for a contour of integration. The
integral $\underset{G^{*}}{\int }f\left( \nu \right) e^{-z\left( t-\nu
\right) }d\nu $, defined over the integration contour $G_{\varepsilon
}^{*}=\left\{ \nu \mid \text{ }\nu \left( \theta \right) =\varepsilon
e^{i\theta }\right\} $ bypassing a singularity of a function $f\left( \nu
\right) e^{-z\left( t-\nu \right) }$ at the point $\nu =0$, is a by-pass
integral. In this case too, similarly to the previous analysis, it is
possible to define an interval (a semi-interval) of a parametric contour
integral convergence of the function $f\left( \nu \right) e^{-z\left( t-\nu
\right) }$. Namely, since in this case 
\begin{equation}
z\left[ t-\nu \left( \theta \right) \right] =z\left[ t-\varepsilon \left(
\cos \theta +i\sin \theta \right) \right] =  \label{25}
\end{equation}
\begin{equation*}
=\func{Re}z\left( t-\varepsilon \cos \theta \right) +\varepsilon \func{Im}%
z\sin \theta +
\end{equation*}
\begin{equation*}
+i\left[ \func{Im}z\left( t-\varepsilon \cos \theta \right) -\varepsilon 
\func{Re}z\sin \theta \right] \text{,}
\end{equation*}
then the condition (\ref{23}), for the contour of integration $G$, reduces
to the condition 
\begin{equation}
\func{Re}z\left[ \left( t-\varepsilon \cos \theta \right) +\varepsilon \frac{%
\func{Im}z}{\func{Re}z}\sin \theta \right] >0\text{\textit{; }}\func{Re}z>0%
\text{,}  \label{26}
\end{equation}
for the contour of integration $G_{\varepsilon }^{*}$.

In view of the fact that $t>\varepsilon >0$, there exists a positive real
number $k$\textit{:} $k\in R_{+}^{1}$, such that $t=\left( 1+k\right)
\varepsilon $. Hence, the condition (\ref{26}) reduces to the condition 
\begin{equation}
\varepsilon \func{Re}z\left[ k+\left( 1-\cos \theta \right) +\tan \varphi
\sin \theta \right] >0\text{\textit{; }}\func{Re}z>0\text{.}  \label{27}
\end{equation}

As $\left| \sin \theta \right| \leq 1$, for $\theta \in \left( -\pi
,0\right) $, the condition (\ref{27}) is satisfied if and only if $\varphi
\in \left( -\frac{\pi }{2},\arctan k\right] $, in other words for each $%
\varphi \in \left( -\frac{\pi }{2},\arctan k\right] $ it holds 
\begin{equation}
\underset{\left| z\right| \rightarrow +\infty }{\lim }z\underset{%
G_{\varepsilon }^{*}}{\int }f\left( \nu \right) e^{-z\left( t-\nu \right)
}d\nu =0\text{.}  \label{28}
\end{equation}

On the other hand, since the real meromorphic function $f\left[ \left( \func{%
Re}\nu \right) \right] $ satisfies in the semi-segment $\left[
t_{0},0\right) $, as well as in the semi-interval $\left( 0,t\right] $,
general well-known \textit{Dirichlet's} conditions, \cite{Mi} and \cite{To},
then, for each $\func{Re}z\geq 0$, it holds\textit{:} 
\begin{equation*}
\underset{\left| z\right| \rightarrow +\infty }{\lim }z\underset{t_{0}}{%
\overset{-\varepsilon }{\int }}f\left( \tau \right) e^{-z\left( t-\tau
\right) }d\tau =0\text{ and}\underset{\left| z\right| \rightarrow +\infty }{%
\lim }z\underset{\varepsilon }{\overset{t}{\int }}f\left( \tau \right)
e^{-z\left( t-\tau \right) }d\tau =f\left( t\right) \text{,}
\end{equation*}
respectively, where $\tau =\func{Re}\nu $. Therefore, the semi-interval $%
\left( -\frac{\pi }{2},\arctan k\right] $ is a semi-interval of convergence
of the parametric contour integral 
\begin{equation}
\underset{\left| z\right| \rightarrow +\infty }{\lim }z\underset{G^{*}}{\int 
}f\left( \nu \right) e^{-z\left( t-\nu \right) }d\nu =f\left( t\right) \text{%
.}  \label{29}
\end{equation}

Taking into consideration the fact that in the limit as $\varepsilon
\rightarrow 0^{+}$ and for each $t>0$\textit{:} $k\rightarrow +\infty $ ( $%
t=\left( 1+k\right) \varepsilon $), the interval $\left( -\frac{\pi }{2},%
\frac{\pi }{2}\right) $ ($\varphi \in \left( -\frac{\pi }{2},\frac{\pi }{2}%
\right) $) becomes an interval of convergence of an improper integral 
\begin{equation}
\underset{\left| z\right| \rightarrow +\infty }{\lim }z\left[ v.t.\underset{%
t_{0}}{\overset{t}{\int }}f\left( \tau \right) e^{-z\left( t-\tau \right)
}d\tau \right] =f\left( t\right)  \label{30}
\end{equation}
absolutely existing in the segment $\left[ t_{0},t\right] $\textit{;} $t>0$.
In other words, for each $\func{Re}z>0$, it holds (\ref{30}).

In the similar manner it can be proved to be 
\begin{equation}
\underset{\left| z\right| \rightarrow +\infty }{\lim }z\left[ v.t.\underset{%
t_{0}}{\overset{t}{\int }}f\left( \tau \right) e^{-z\left( \tau
-t_{0}\right) }d\tau \right] =f\left( t_{0}\right) \text{,}  \label{31}
\end{equation}
for each $\func{Re}z>0$ and $t>0$ ($t\in \left( t_{0},t_{1}\right) $.

\begin{description}
\item[Comment]  It should be emphasized that as distinguished from limiting
values of contour integrals\textit{:} 
\begin{equation*}
\underset{\left| z\right| \rightarrow +\infty }{\lim }z\underset{G}{\int }%
f\left( \nu \right) e^{-z\left( t-\nu \right) }d\nu \text{ and}\underset{%
\left| z\right| \rightarrow +\infty }{\lim }z\underset{G^{*}}{\int }f\left(
\nu \right) e^{-z\left( t-\nu \right) }d\nu \text{,}
\end{equation*}
which are equal, their intervals of convergence are different.$%
\blacktriangledown $
\end{description}

\subsubsection{\textit{Fourier} formula\label{Fourier}}

Based on the results obtained by previous analysis one may say that in
addition to a class of real valued functions satisfying so-called general%
\textit{\ Dirichlet's} conditions in the real axis segment $\left[
t_{0},t_{1}\right] $, for which the functional expressions\textit{:} (\ref
{12}) and (\ref{13}), converge for each $\func{Re}z\geq 0$, \cite{M-K (78)},
there exists an one other class of real valued functions for which \textit{%
Cauchy's} formula is still in effect, a class of the real valued meromorphic
functions whose finitely many isolated singularities lie onto the segment $%
\left[ t_{0},t_{1}\right] $.

From \textit{Cauchy's} formula, and for $p\left( z\right) =-1$ and $w\left(
z\right) =e^{az}$ ($a\in R_{+}^{1}$) i.e. $q\left( z\right) =e^{az}-1$, it
is immediately obtained that for $t_{1}-a<t<t_{0}+a$ and $t\neq t_{si}$%
\begin{equation}
f\left( t\right) =\frac{1}{a}\left[ v.t.\overset{t_{1}}{\underset{t_{0}}{%
\int }}f\left( \tau \right) d\tau \right] +\frac{2}{a}\overset{+\infty }{%
\underset{k=1}{\sum }}\left[ v.t.\overset{t_{1}}{\underset{t_{0}}{\int }}%
f\left( \tau \right) \cos \frac{2k\pi \left( t-\tau \right) }{a}d\tau
\right] \text{.}  \label{32}
\end{equation}

If $a=2\pi $, $t_{0}=-\pi $ and $t_{1}=\pi $, then the equation (\ref{32})
represents an expansion of a real valued meromorphic function $t\mapsto
f\left( t\right) $ into a \textit{Fourier} trigonometric series in the
interval $\left( -\pi ,\pi \right) $, more exactly, for each $-\pi <t<\pi $
and $t\neq t_{si}$, where $t_{si}$ are break points of the function $%
t\mapsto f\left( t\right) $ in the interval $\left( -\pi ,\pi \right) $, it
holds 
\begin{equation}
f\left( t\right) =\frac{1}{2}A_{0}+\overset{+\infty }{\underset{k=1}{\sum }}%
\left[ A_{k}\cos \left( kt\right) +B_{k}\sin \left( kt\right) \right] \text{,%
}  \label{33}
\end{equation}
where 
\begin{equation}
A_{k}=\frac{1}{\pi }\left[ v.t.\overset{\pi }{\underset{-\pi }{\int }}%
f\left( \tau \right) \cos \left( k\tau \right) d\tau \right] \text{\textit{;}
}k\in N_{0}  \label{34}
\end{equation}
and 
\begin{equation}
B_{k}=\frac{1}{\pi }\left[ v.t.\overset{\pi }{\underset{-\pi }{\int }}%
f\left( \tau \right) \sin \left( k\tau \right) d\tau \right] \text{\textit{;}
}k\in N\text{.}  \label{35}
\end{equation}

\begin{description}
\item[Comment]  According to the well-known result of \textit{Dirichlet's}
theorem, see - \textit{Theorem 1, Section 5.4, Chapter 5, p. 65}, \cite{Mi}
- in the general case of a class of real valued functions $t\mapsto
f_{d}\left( t\right) $ satisfying the general \textit{Dirichlet's}
conditions in the segment $\left[ -\pi ,\pi \right] $, the \textit{Fourier}
trigonometric series on the right-hand side of the equation (\ref{33}) can
be said to converge to a function $F_{d}\left( t\right) $. Clearly, at all
points of the interval $\left( -\pi ,\pi \right) $ at which the function $%
t\mapsto f_{d}\left( t\right) $ is continuous, a convergent value of \textit{%
Fourier} series is equal to the function value: $F_{d}\left( t\right)
=f_{d}\left( t\right) $. Considering the consequence of whether \textit{%
Bessel's} inequality \cite{Mi} or \textit{Riemann-Lebesque's} theorem - 
\textit{Theorem 2, Section 6.2, Chapter 5, p. 96,} \cite{Mi} - \textit{%
Fourier's} coefficients of the function $t\mapsto f_{d}\left( t\right) $%
\textit{:} $A_{k}$ and $B_{k}$, tend to zero as $k\rightarrow +\infty $.
This is important from the viewpoint of the convergence of infinite
numerical series obtained by expansion of functions $t\mapsto f_{d}\left(
t\right) $ into \textit{Fourier} trigonometric series. Note that the
conditions of \textit{Dirichlet's} theorem are only sufficient conditions
for convergence of \textit{Fourier} trigonometric series of functions $%
t\mapsto f_{d}\left( t\right) $.

A nature of \textit{Fourier} trigonometric series convergence of a class of
real valued meromorphic functions $t\mapsto f_{m}\left( t\right) $, whose
finitely many isolated singularities lie onto the segment $\left[ -\pi ,\pi
\right] $, can be said to be different from a case to a case. The same holds
also for \textit{Fourier's} coefficients of a function $t\mapsto f_{m}\left(
t\right) $ in the limit as $k\rightarrow +\infty $. Namely, in the general
case of real meromorphic functions $t\mapsto f_{m}\left( t\right) $, at all
points of the interval $\left( -\pi ,\pi \right) $ at which a function $%
t\mapsto f_{m}\left( t\right) $ is continuous, \textit{Fourier}
trigonometric series of the function $t\mapsto f_{m}\left( t\right) $ is
summable, more exactly it has defined sum $F_{m}\left( t\right) $\textit{:} $%
F_{m}\left( t\right) =f_{m}\left( t\right) $. The concept of the sum of 
\textit{Fourier} trigonometric series, in this case, is generalization of
the concept of the sum in the \textit{Caushy's} sense. At the break points $%
t_{si}$ of the function $t\mapsto f_{m}\left( t\right) $ in the segment $%
\left[ -\pi ,\pi \right] $, the sum of \textit{Fourier} trigonometric series
of the real valued meromorphic function $t\mapsto f_{m}\left( t\right) $, in
the general case, is not defined.$\blacktriangledown $
\end{description}

\section{Examples}

\subsection{ Example 1}

\textit{An expansion of the function }$t\mapsto \frac{1}{2}\frac{\sin t}{%
1-\cos t}$\textit{\ into a Fourier trigonometric series in the segment }$%
\left[ -\pi ,\pi \right] $. The function $f\left( t\right) =\frac{1}{2}\frac{%
\sin t}{1-\cos t}$ having at the point a simple pole as a singularity is a
real valued meromorphic function in the segment $\left[ -\pi ,\pi \right] $. 
\textit{Caushy's} principal value ($v.p.$) of an improper integral of the
function $f\left( t\right) $ is equal to\textit{:} 
\begin{equation}
v.p.\overset{\pi }{\underset{-\pi }{\int }}f\left( t\right) dt=\frac{1}{2}%
\underset{\varepsilon \rightarrow 0^{+}}{\lim }\left[ \overset{-\varepsilon 
}{\underset{-\pi }{\int }}\frac{\sin t}{1-\cos t}dt+\overset{\pi }{\underset{%
\varepsilon }{\int }}\frac{\sin t}{1-\cos t}dt\right] =0\text{.}  \label{36}
\end{equation}

As $\underset{z\rightarrow 0}{\lim }\frac{1}{2}\frac{z\sin z}{1-\cos z}=1$%
\textit{;} $z\in C^{1}$, the by-pass integral value in the limit as $%
\varepsilon \rightarrow 0^{+}$ (\textit{Jordan's} singular value \textit{(}$%
v.s.$\textit{)} of the improper integral) is equal to\textit{:} 
\begin{equation}
\underset{\varepsilon \rightarrow 0^{+}}{\lim }\underset{G_{\varepsilon
_{\kappa }}}{\int }\frac{1}{2}\frac{\sin z}{1-\cos z}dz=\left\{ 
\begin{array}{l}
-i\pi \text{\textit{;} }\kappa =1 \\ 
i\pi \text{\textit{;} }\kappa =2
\end{array}
\right. \text{,}  \label{37}
\end{equation}
dependently on the choice of a circular arc $G_{\varepsilon _{\kappa }}$
bypassing the singularity of the function $f\left( t\right) $ in the complex
plane\textit{:} $G_{\varepsilon _{\kappa }}=\left\{ z\mid z=\varepsilon
e^{i\theta _{\kappa }}\text{\textit{;} }\theta _{\kappa }\in \left\{ 
\begin{array}{l}
\left[ -\pi ,0\right] \text{\textit{;} }\kappa =1 \\ 
\left[ \pi ,0\right] \text{\textit{;} }\kappa =2
\end{array}
\right. \right\} $\textit{.}

The total value ($v.t.$) of the improper integral, as a sum of \textit{%
Caushy's} principal value ($v.p.$) and \textit{Jordan's} singular value ($%
v.s.$)\textit{:} 
\begin{equation}
\frac{1}{\pi }v.t.\overset{\pi }{\underset{-\pi }{\int }}f\left( t\right) dt=%
\frac{1}{2\pi }v.t.\overset{\pi }{\underset{-\pi }{\int }}\frac{\sin t}{%
1-\cos t}dt=\left\{ 
\begin{array}{l}
-i \\ 
i
\end{array}
\right.  \label{38}
\end{equation}
absolutely exists in this case and as one can see is not unique.

On the other hand, since $\frac{1}{\pi }\overset{\pi }{\underset{0}{\int }}%
\frac{\sin \left[ \left( k+\frac{1}{2}\right) t\right] }{\sin \left( \frac{1%
}{2}t\right) }dt=1$\textit{;} $k\in N$ - \textit{Formula (6) Section 6.2,
Chapter 6, p. 95}, \cite{Mi} - and $\frac{1}{2}\overset{\pi }{\underset{-\pi 
}{\int }}\cos \left( kt\right) dt=0$\textit{;} $k\in N$, as well as $%
\underset{z\rightarrow 0}{\lim }z\frac{\sin z\cos \left( kz\right) }{2\left(
1-\cos z\right) }=1$, then it follows that 
\begin{equation}
B_{k}=\frac{1}{\pi }v.t.\overset{\pi }{\underset{-\pi }{\int }}\frac{\sin
t\sin \left( kt\right) }{2\left( 1-\cos t\right) }dt=\frac{1}{2\pi }v.t.%
\overset{\pi }{\underset{-\pi }{\int }}\cot \frac{t}{2}\sin \left( kt\right)
dt=1\text{\textit{; }}k\in N\text{,}  \label{39}
\end{equation}
as well as 
\begin{equation}
A_{k}=\frac{1}{\pi }v.t.\overset{\pi }{\underset{-\pi }{\int }}\frac{\sin
t\cos \left( kt\right) }{2\left( 1-\cos t\right) }dt=\left\{ 
\begin{array}{l}
-i \\ 
i
\end{array}
\right. \text{\textit{;} }k\in N\text{,}  \label{40}
\end{equation}

Accordingly, and by the \textit{Fourier} formula (\ref{33}), a \textit{%
Fourier} trigonometric series of the function $t\mapsto \frac{1}{2}\frac{%
\sin t}{1-\cos t}$ can be expressed by the following functional form 
\begin{equation}
\frac{1}{2}\frac{\sin t}{1-\cos t}=\underset{k=1}{\overset{+\infty }{\sum }}%
\sin \left( kt\right) \mp \frac{i}{2}\left[ 1+2\underset{k=1}{\overset{%
+\infty }{\sum }}\cos \left( kt\right) \right] \text{,}  \label{41}
\end{equation}
that is, the equalities 
\begin{equation}
\frac{1}{2}\frac{\sin t}{1-\cos t}=\underset{k=1}{\overset{+\infty }{\sum }}%
\sin \left( kt\right) \text{ and }1+2\underset{k=1}{\overset{+\infty }{\sum }%
}\cos \left( kt\right) =0\text{ i.e.}  \label{42}
\end{equation}
\begin{equation*}
\frac{1}{2}\frac{1-e^{\pm it}}{1-\cos t}=-\underset{k=1}{\overset{+\infty }{%
\sum }}e^{\pm ikt}
\end{equation*}
hold for each $t\in \left( -\pi ,\pi \right) $ and $t\neq 0$, respectively.

By relation (\ref{15}), for $t=\pm \pi $, it follows that 
\begin{equation}
\underset{k=1}{\overset{+\infty }{\sum }}\sin \left( k\pi \right) =0\text{
and }1+2\underset{k=1}{\overset{+\infty }{\sum }}\cos \left( k\pi \right)
=1+2\underset{k=1}{\overset{+\infty }{\sum }}\left( -1\right) ^{k}=0\text{.}%
\blacktriangledown  \label{43}
\end{equation}

\subsection{ Example 2}

\textit{An expansion of the function }$t\mapsto \frac{1}{2}\frac{1}{1-\cos t}
$\textit{\ into a Fourier trigonometric series in the segment }$\left[ -\pi
,\pi \right] $. The real valued meromorphic function $f\left( t\right) =%
\frac{1}{2}\frac{1}{1-\cos t}$ in the segment $\left[ -\pi ,\pi \right] $
has the second order pole at the point $t=0$ as a singularity. The improper
integral $\frac{1}{2}\underset{-\pi }{\overset{\pi }{\int }}\frac{dt}{1-\cos
t}$ absolutely exists and reduces to the indefinite expression of difference
of infinities $\infty -\infty $. Namely, independently on the choice of the
circular arc bypassing the singularity $z=0$ of the function $z\mapsto \frac{%
1}{2}\frac{1}{1-\cos z}$ in the complex plane (for example $G_{\varepsilon
}=\left\{ z\mid z=\varepsilon e^{i\theta }\text{\textit{;} }\theta \in
\left[ -\pi ,0\right] \right\} $), on the one hand it holds 
\begin{equation}
\frac{1}{2}v.t.\underset{-\pi }{\overset{\pi }{\int }}\frac{dt}{1-\cos t}=%
\frac{1}{2}\underset{\varepsilon \rightarrow 0^{+}}{\lim }\left[ \overset{%
-\varepsilon }{\underset{-\pi }{\int }}\frac{dt}{1-\cos t}+\underset{%
G_{\varepsilon }}{\int }\frac{dz}{1-\cos z}+\overset{\pi }{\underset{%
\varepsilon }{\int }}\frac{dt}{1-\cos t}\right] =  \label{44}
\end{equation}
\begin{equation*}
=\underset{\varepsilon \rightarrow 0^{+}}{\lim }\left[ \frac{\sin
\varepsilon }{1-\cos \varepsilon }+\frac{1}{2}\underset{G_{\varepsilon }}{%
\int }\frac{dz}{1-\cos z}\right] \text{.}
\end{equation*}

On the other hand, since - \textit{see} \textit{Definition 4, Section 2.2,
Chapter 2, p. 82}, \cite{St} 
\begin{equation}
\underset{G_{\varepsilon }}{\int }\frac{dz}{1-\cos z}=\overset{0}{\underset{%
-\pi }{\int }}\frac{\frac{d}{d\theta }\left[ z\left( \theta \right) \right] 
}{1-\cos \left[ z\left( \theta \right) \right] }d\theta =\overset{0}{%
\underset{-\pi }{\int }}\frac{i\varepsilon e^{i\theta }}{1-\cos \left(
\varepsilon e^{i\theta }\right) }d\theta \text{,}  \label{45}
\end{equation}
that is 
\begin{equation}
\underset{G_{\varepsilon }}{\int }\frac{dz}{1-\cos z}=-\frac{\sin \left[
z\left( \theta \right) \right] }{1-\cos \left[ z\left( \theta \right)
\right] }\left| _{-\pi }^{0}\right. =-\frac{\sin \left( \varepsilon
e^{i\theta }\right) }{1-\cos \left( \varepsilon e^{i\theta }\right) }\left|
_{-\pi }^{0}\right. =-\frac{2\sin \varepsilon }{1-\cos \varepsilon }\text{,}
\label{46}
\end{equation}
then finally it follows that the total value ($v.t.$) of the improper
integral $\underset{-\pi }{\overset{\pi }{\int }}\frac{dt}{2\left( 1-\cos
t\right) }$ is equal to the value zero\textit{:} 
\begin{equation}
v.t.\underset{-\pi }{\overset{\pi }{\int }}\frac{dt}{2\left( 1-\cos t\right) 
}=0\text{.}  \label{47}
\end{equation}

The complex function $z\mapsto \frac{z^{k}}{z-1}$\textit{;} $k\in N$ of
complex variable $z$ is a meromorphic function having at the point $z=1$ a
simple pole as singularity. Since $\underset{z\rightarrow 1}{\lim }\left[
\left( z-1\right) \frac{z^{k}}{z-1}\right] =1$ then, according to the result
(\ref{1}) in the \textit{Section 2.1} of the paper, \textit{Cauchy's}
principle value ($v.p.$) of the improper integral $\underset{G}{\int }\frac{%
z^{k}}{z-1}dz$ over the circular contour of integration $G$\textit{:} $%
G=\left\{ z\mid z=\varepsilon e^{i\theta }\text{\textit{;} }\theta \in
\left[ -\pi ,\pi \right] \right\} $ is equal to\textit{:} $v.p.\underset{G}{%
\int }\frac{z^{k}}{z-1}dz=i\pi $.

With regard to the fact that $z=\varepsilon e^{i\theta }$ onto the
integration contour $G$, it follows that 
\begin{equation}
v.p.\underset{G}{\int }\frac{z^{k}}{z-1}dz=v.p.\underset{-\pi }{\overset{\pi 
}{\int }}\frac{ie^{ik\theta }}{1-e^{-i\theta }}d\theta =v.p.\left[ \underset{%
-\pi }{\overset{\pi }{\int }}\frac{i\cos \left( k\theta \right) }{2\left(
1-\cos \theta \right) }d\theta -\right.  \label{48}
\end{equation}
\begin{equation*}
\left. -\underset{-\pi }{\overset{\pi }{\int }}\frac{i\cos \left[ \left(
k+1\right) \theta \right] }{2\left( 1-\cos \theta \right) }d\theta \right]
=i\pi \text{,}
\end{equation*}
since $v.p.\underset{-\pi }{\overset{\pi }{\int }}\frac{\cos \left( k\theta
\right) \sin \theta }{1-\cos \theta }d\theta =0$ and $\underset{-\pi }{%
\overset{\pi }{\int }}\sin \left( k\theta \right) d\theta =0$\textit{;} for
each $k\in N$.

On the other hand, for $k\in N$\textit{:} $\underset{z\rightarrow 0}{\lim }%
\left[ z\cos \left( kz\right) \right] =0$ and $\underset{z\rightarrow 0}{%
\lim }\frac{z\sin \left( kz\right) \sin z}{1-\cos z}=0$. According to the
result of \textit{Jordan's} lemma - \textit{Theorem 1, Subsection 3.1.4,
Section 3.1, Chapter 3, p. 52}, \cite{M-K (78)} - it holds 
\begin{equation}
\frac{1}{2\pi i}\underset{\varepsilon \rightarrow 0^{+}}{\lim }\underset{%
G_{\varepsilon _{\kappa }}}{\int }\left\{ \frac{\cos \left( kz\right) -\cos
\left[ \left( k+1\right) z\right] }{2\left( 1-\cos z\right) }\right\} dz=
\label{49}
\end{equation}
\begin{equation*}
=\frac{1}{4\pi i}\underset{\varepsilon \rightarrow 0^{+}}{\lim }\underset{%
G_{\varepsilon _{\kappa }}}{\int }\left[ \cos \left( kz\right) -\frac{\sin
\left( kz\right) \sin z}{\left( 1-\cos z\right) }\right] dz=0\text{,}
\end{equation*}
where the singularity $z=0$ of the meromorphic function $z\mapsto \frac{1}{2}%
\frac{1}{1-\cos z}$ is bypassed by the parts $G_{\varepsilon _{\kappa }}$%
\textit{:} $G_{\varepsilon _{\kappa }}=\left\{ z\mid z=e^{i\theta _{\kappa }}%
\text{\textit{;} }\theta _{\kappa }\in \left\{ 
\begin{array}{l}
\left[ \pi ,0\right] \text{\textit{;} }\kappa =1 \\ 
\left[ -\pi ,0\right] \text{\textit{;} }\kappa =2
\end{array}
\right. \right\} $, of a circular path of integration\textit{:} $%
G_{\varepsilon }=\left\{ z\mid z=e^{i\theta }\text{\textit{;} }\theta \in
\left[ -\pi ,\pi \right] \right\} $.

Finally, from results\textit{:} (\ref{48}) and (\ref{49}), the integral
equation is obtained 
\begin{equation}
\frac{1}{2\pi }v.t.\left[ \underset{-\pi }{\overset{\pi }{\int }}\frac{\cos
\left( k\theta \right) }{2\left( 1-\cos \theta \right) }d\theta -\underset{%
-\pi }{\overset{\pi }{\int }}\frac{\cos \left[ \left( k+1\right) \theta
\right] }{2\left( 1-\cos \theta \right) }d\theta \right] =\frac{1}{2}\text{%
\textit{;} }k\in N\text{,}  \label{50}
\end{equation}
and that is in agreement with result (\ref{39}).

Further, since $\cos \left( 2\theta \right) =1-2\left( \sin \theta \right)
^{2}$, it follows that 
\begin{equation}
v.t.\underset{-\pi }{\overset{\pi }{\int }}\frac{\cos \left( 2\theta \right) 
}{2\left( 1-\cos \theta \right) }d\theta =v.t.\left[ \underset{-\pi }{%
\overset{\pi }{\int }}\frac{d\theta }{2\left( 1-\cos \theta \right) }-%
\underset{-\pi }{\overset{\pi }{\int }}\frac{\left( \sin \theta \right) ^{2}%
}{1-\cos \theta }d\theta \right] \text{,}  \label{51}
\end{equation}
that is $\frac{1}{\pi }v.t.\underset{-\pi }{\overset{\pi }{\int }}\frac{\cos
\left( 2\theta \right) }{2\left( 1-\cos \theta \right) }d\theta =-2$, in
view of the results\textit{:} (\ref{39}) and (\ref{47}).

Consequently, considering (\ref{50}) it has been just proved by a method of
mathematical induction that for each $k\in N$ it holds 
\begin{equation}
A_{k}=\frac{1}{\pi }v.t.\underset{-\pi }{\overset{\pi }{\int }}\frac{\cos
\left( kt\right) }{2\left( 1-\cos t\right) }dt=-k\text{.}  \label{52}
\end{equation}

As for an improper integral $\underset{-\pi }{\overset{\pi }{\int }}\frac{%
\sin \left( kt\right) }{2\left( 1-\cos t\right) }dt$, its total value ($v.t.$%
) reduces to the value of by-pass integrals $\underset{G_{\varepsilon
_{\kappa }}}{\int }\frac{\sin \left( kz\right) }{2\left( 1-\cos z\right) }dz$
in the limit as $\varepsilon \rightarrow 0^{+}$. Taking into account that
for each $k\in N$\textit{:} $\underset{z\rightarrow 0}{\lim }\frac{z\sin
\left( kz\right) }{2\left( 1-\cos z\right) }=k$, it follows that 
\begin{equation}
B_{k}=\frac{1}{\pi }v.t.\underset{-\pi }{\overset{\pi }{\int }}\frac{\sin
\left( kt\right) }{2\left( 1-\cos t\right) }dt=  \label{53}
\end{equation}
\begin{equation*}
=\frac{1}{\pi }\underset{\varepsilon \rightarrow 0^{+}}{\lim }\underset{%
G_{\varepsilon _{\kappa }}}{\int }\frac{\sin \left( kz\right) }{2\left(
1-\cos z\right) }dz=\left\{ 
\begin{array}{l}
-ik\text{\textit{;} }\kappa =1 \\ 
ik\text{\textit{;} }\kappa =2
\end{array}
\right. \text{\textit{; }}k\in N\text{. }
\end{equation*}

Therefore, \textit{Fourier} trigonometric series of the function $t\mapsto 
\frac{1}{2}\frac{1}{1-\cos t}$ in the segment $\left[ -\pi ,\pi \right] $,
according to results: (\ref{47}) and (\ref{52}) as well as (\ref{53}), can
be expressed by the following functional form 
\begin{equation}
\frac{1}{2}\frac{1}{1-\cos t}=-\overset{+\infty }{\underset{k=1}{\sum }}%
k\cos \left( kt\right) \mp i\overset{+\infty }{\underset{k=1}{\sum }}k\sin
\left( kt\right) =-\overset{+\infty }{\underset{k=1}{\sum }}ke^{\pm ikt}%
\text{,}  \label{54}
\end{equation}
that is, the equalities 
\begin{equation}
\frac{1}{2}\frac{1}{1-\cos t}=-\overset{+\infty }{\underset{k=1}{\sum }}%
k\cos \left( kt\right) \text{ and }\overset{+\infty }{\underset{k=1}{\sum }}%
k\sin \left( kt\right) =0\text{,}  \label{55}
\end{equation}
hold for each $t\in \left( -\pi ,\pi \right) $ and $t\neq 0$, respectively.

In the extreme points of the segment $\left[ -\pi ,\pi \right] $, from (\ref
{15}), it follows that 
\begin{equation}
\overset{+\infty }{\underset{k=1}{\sum }}k\cos \left( k\pi \right) =\overset{%
+\infty }{\underset{k=1}{\sum }}k\left( -1\right) ^{k}=-\frac{1}{4}\text{
and }\overset{+\infty }{\underset{k=1}{\sum }}k\sin \left( k\pi \right) =0%
\text{.}\blacktriangledown  \label{56}
\end{equation}

\begin{description}
\item[Comment]  By an expansion of the real valued functions of real
variable $t$\textit{:} 
\begin{equation*}
f\left( t\right) =\left\{ 
\begin{array}{l}
\frac{\sin t}{2\left( 1-\cos t\right) }\text{\textit{;} }\tau _{0}\leq
\left| t\right| \leq \pi  \\ 
0\text{\textit{;} }\left| t\right| <\tau _{0}
\end{array}
\right. \text{ and }g\left( t\right) =\left\{ 
\begin{array}{l}
b\text{\textit{;} }\tau _{0}\leq t\leq \pi  \\ 
0\text{\textit{;} }\left| t\right| <\tau _{0} \\ 
a\text{\textit{;} }-\pi \leq t\leq -\tau _{0}
\end{array}
\right. \text{\textit{;} }\tau _{0}>0\text{,}
\end{equation*}
satisfying \textit{Dirichlet's} conditions in the segment $\left[ -\pi ,\pi
\right] $, into \textit{Fourier} trigonometric series, it is obtained that 
\begin{equation*}
f\left( t\right) =\overset{+\infty }{\underset{k=1}{\sum }}\left[ \frac{1}{%
\pi }\underset{\tau _{0}}{\overset{\pi }{\int }}\frac{\sin \tau }{1-\cos
\tau }\sin \left( k\tau \right) d\tau \right] \sin \left( kt\right) \text{%
\textit{; }}\left| t\right| \in \left( \tau _{0},\pi \right) \text{,}
\end{equation*}
\begin{equation*}
g\left( t\right) =\frac{1}{2\pi }\left( \underset{-\pi }{\overset{-\tau _{0}%
}{\int }}ad\tau +\underset{\tau _{0}}{\overset{\pi }{\int }}bd\tau \right) +%
\overset{+\infty }{\underset{k=1}{\sum }}\frac{1}{\pi }\left\{ \left[ 
\underset{-\pi }{\overset{-\tau _{0}}{\int }}a\sin \left( k\tau \right)
d\tau +\underset{\tau _{0}}{\overset{\pi }{\int }}b\sin \left( k\tau \right)
d\tau \right] \sin \left( kt\right) +\right. 
\end{equation*}
\begin{equation*}
\left. +\left[ \underset{-\pi }{\overset{-\tau _{0}}{\int }}a\cos \left(
k\tau \right) d\tau +\underset{\tau _{0}}{\overset{\pi }{\int }}b\cos \left(
k\tau \right) d\tau \right] \cos \left( kt\right) \right\} \text{\textit{;} }%
\left| t\right| \in \left( \tau _{0},\pi \right) \text{,}
\end{equation*}
that is\footnote{{\footnotesize By the well-known trigonometrical equalities%
\textit{:} }$\sin \left[ \left( k+1\right) t\right] =\sin \left( kt\right)
\cos t+\cos \left( kt\right) \sin t${\footnotesize \ and }$\frac{\sin t\sin
\left[ \left( k+1\right) t\right] }{1-\cos t}=\frac{\sin t\cos t}{1-\cos t}%
\sin \left( kt\right) +\left( 1+\cos t\right) \cos \left( kt\right) $%
{\footnotesize , as well as }$\sin \left( kt\right) \cos t=\frac{1}{2}%
\left\{ \sin \left[ \left( k-1\right) t\right] +\sin \left[ \left(
k+1\right) t\right] \right\} ${\footnotesize \ and }$\cos t\cos \left(
kt\right) =\frac{1}{2}\left\{ \cos \left[ \left( k-1\right) t\right] +\cos
\left[ \left( k+1\right) t\right] \right\} ${\footnotesize , the following
recurrent formula for \textit{Fourier's} coefficients of the function }$%
f\left( t\right) ${\footnotesize \ in the segment }$\left[ -\pi ,\pi \right] 
$ {\footnotesize is obtained } 
\begin{equation*}
\frac{1}{\pi }\underset{\tau _{0}}{\overset{\pi }{\int }}\frac{\sin t\sin
\left[ \left( k+1\right) t\right] }{1-\cos t}dt=\frac{1}{\pi }\underset{\tau
_{0}}{\overset{\pi }{\int }}\frac{\sin t\sin \left[ \left( k-1\right)
t\right] }{1-\cos t}dt-\frac{2\sin \left( k\tau _{0}\right) }{k\pi }-
\end{equation*}
\begin{equation*}
-\frac{\sin \left[ \left( k+1\right) \tau _{0}\right] }{\left( k+1\right)
\pi }-\frac{\sin \left[ \left( k-1\right) \tau _{0}\right] }{\left(
k-1\right) \pi }\text{\textit{; }}k\in N
\end{equation*}
\par
{\footnotesize On the other hand, for }$k=1${\footnotesize , that is, for }$%
k=2${\footnotesize , it holds\textit{:} }$\frac{1}{\pi }\underset{\tau _{0}}{%
\overset{\pi }{\int }}\frac{\left( \sin t\right) ^{2}}{1-\cos t}dt=\frac{1}{%
\pi }\underset{\tau _{0}}{\overset{\pi }{\int }}\left( 1+\cos t\right) dt=1-%
\frac{\tau _{0}}{\pi }-\frac{\sin \tau _{0}}{\pi }${\footnotesize , that is, 
}$\frac{1}{\pi }\underset{\tau _{0}}{\overset{\pi }{\int }}\frac{\sin \left(
2t\right) \sin t}{1-\cos t}dt=\frac{2}{\pi }\underset{\tau _{0}}{\overset{%
\pi }{\int }}\left( 1+\cos t\right) \cos tdt=1-\frac{\tau _{0}}{\pi }-2\frac{%
\sin \tau _{0}}{\pi }${\footnotesize \ }$-\frac{\sin \left( 2\tau
_{0}\right) }{2\pi }${\footnotesize , respectively.}} 
\begin{equation}
f\left( t\right) =\overset{+\infty }{\underset{k=1}{\sum }}\left[ 1+\frac{%
\tau _{0}}{\pi }-2\overset{k-1}{\underset{\kappa =0}{\sum }}\frac{\sin
\left( \kappa \tau _{0}\right) }{\kappa \pi }-\frac{\sin \left( k\tau
_{0}\right) }{k\pi }\right] \sin \left( kt\right) \text{\textit{;} }\left|
t\right| \in \left( \tau _{0},\pi \right) \text{,}  \label{57}
\end{equation}
\begin{equation}
g\left( t\right) =\frac{a+b}{2}-\frac{a+b}{2\pi }\left[ \frac{1}{2}+\overset{%
+\infty }{\underset{k=1}{\sum }}\frac{\sin \left( k\tau _{0}\right) }{k\tau
_{0}}\cos \left( kt\right) \right] \tau _{0}+  \label{58}
\end{equation}
\begin{equation*}
+\frac{b-a}{\pi }\overset{+\infty }{\underset{k=1}{\sum }}\left[ \cos \left(
k\tau _{0}\right) -\left( -1\right) ^{k}\right] \frac{\sin \left( kt\right) 
}{k}\text{\textit{;} }\left| t\right| \in \left( \tau _{0},\pi \right) .
\end{equation*}

From the functional relation (\ref{58}) it follows for $a=b$, $\left|
t\right| \in \left( \tau _{0},\pi \right) $ and $\tau _{0}>0$ that 
\begin{equation}
\frac{1}{2}+\overset{+\infty }{\underset{k=1}{\sum }}\frac{\sin \left( k\tau
_{0}\right) }{k\tau _{0}}\cos \left( kt\right) =0\text{.}  \label{59}
\end{equation}

Thus, for $t\in \left( \tau _{0},\pi \right) $ and $\tau _{0}>0$ it holds 
\begin{equation}
\frac{\pi }{2}=\overset{+\infty }{\underset{k=1}{\sum }}\left[ \cos \left(
k\tau _{0}\right) -\left( -1\right) ^{k}\right] \frac{\sin \left( kt\right) 
}{k}\text{,}  \label{60}
\end{equation}
that is 
\begin{equation}
\overset{+\infty }{\underset{k=1}{\sum }}\cos \left( k\tau _{0}\right) \frac{%
\sin \left( kt\right) }{k}=\frac{\pi }{2}-\frac{t}{2}\text{\textit{;} }t\in
\left( \tau _{0},\pi \right) \text{,}  \label{61}
\end{equation}
since for $t\in \left( -\pi ,\pi \right) $, \cite{Mi}, 
\begin{equation}
\overset{+\infty }{\underset{k=1}{\sum }}\left( -1\right) ^{k}\frac{\sin
\left( kt\right) }{k}=-\frac{t}{2}\text{.}  \label{62}
\end{equation}

On the other hand, taking into account the fact that 
\begin{equation*}
\underset{k\rightarrow +\infty }{\lim }\left[ 1+\frac{\tau _{0}}{\pi }-2%
\overset{k-1}{\underset{\kappa =0}{\sum }}\frac{\sin \left( \kappa \tau
_{0}\right) }{\kappa \pi }-\frac{\sin \left( k\tau _{0}\right) }{k\pi }%
\right] =0\text{,}
\end{equation*}
see the last \textit{Comment} in the preceding \textit{Section} of this
paper, finally it follows for $\tau _{0}\in \left( 0,\pi \right) $ that 
\begin{equation}
\overset{+\infty }{\underset{k=0}{\sum }}\frac{\sin \left( k\tau _{0}\right) 
}{k}=\frac{\pi }{2}+\frac{\tau _{0}}{2}\text{.}  \label{63}
\end{equation}

Since a real parameter $\tau _{0}$ takes any value from the interval $\left(
0,\pi \right) $, even if that is finitely small value, it would be
reasonable to ask\textit{:} Whether the functional expressions\textit{:} (%
\ref{57}) and (\ref{59}) as well as (\ref{60}) and (\ref{61}), hold in the
limit as $\tau _{0}\rightarrow 0^{+}$? In other words\textit{:} Are the
limiting values of sums, in these emphasized cases, equal to sums of
limiting values of the functional expressions, as $\tau _{0}\rightarrow 0^{+}
$, respectively? On the basis of previously derived results in the \textit{%
Example 1} and of the above obtained result (\ref{63}) as well as of the
well-known result of the series theory\textit{:} $\frac{\pi }{4}=\overset{%
+\infty }{\underset{k=1}{\sum }}\frac{\sin \left[ \left( 2k-1\right)
t\right] }{2k-1}$\textit{;} for $t\in \left( 0,\pi \right) $, \cite{Mi}, an
answer to the former questions is yes. However, the problem of
generalization of a preceding conclusion stays open and can be a subject of
a separate analysis.

Similarly, since $\frac{d}{dt}\left[ \frac{\sin t}{2\left( 1-\cos t\right) }%
\right] =-\frac{1}{2\left( 1-\cos t\right) }$ and $\frac{d}{dt}\left\{ \frac{%
1}{2}\ln \left[ 2\left( 1-\cos t\right) \right] \right\} =\frac{\sin t}{%
2\left( 1-\cos t\right) }$ for $\left| t\right| \in \left( 0,\pi \right) $,
then closely related to results\textit{:} (\ref{42}) and (\ref{55}), of the
paper, as well as to the well-known result of the series theory\textit{:} $-%
\frac{1}{2}\ln \left[ 2\left( 1-\cos t\right) \right] =\overset{+\infty }{%
\underset{k=1}{\sum }}\frac{\cos \left( kt\right) }{k}$, for $\left|
t\right| \in \left( 0,\pi \right) $, \cite{Sl}, is the following question%
\textit{:} In which general cases the derivative of a sum of infinite
functional series is equal to the sum of the derivative of any series
member, separately?

This question also stays open for a separate analysis.$\blacktriangledown $
\end{description}

\section{Conclusion}

Taking into consideration the fact that obtained results are theoretical
news, one can say that the certain possibilities for expansion of some
mathematical analysis knowledge connecting to the problems to which a proper
attention has been paid in this paper are opening up. Thus, from viewpoint
of results, derived in the \textit{Subsubsection \ref{Fourier}} of the paper
for instance, and having in mind the fact that causality related to the area
of \textit{Fourier} trigonometric series of real valued functions is the
theory of partial differential equations, it is obvious which possibilities
are opening up in this area of mathematics.

On the other hand, disregarding the fact that the results of the paper are,
in a certain sense, the theoretical news, some of them have been
predictable. So, the alternative numerical series\textit{:} $\overset{%
+\infty }{\underset{k=0}{\sum }}\left( -1\right) ^{k}$, has the defined sum,
more exactly it is summable and its sum is equal to $\frac{1}{2}$, just as
it has been assumed yet by \textit{Euler} and \textit{Dalamber}. Making use
of this assumption they obtained absolutely exact results. It is nothing
other to be left than to prove validity of this assumption. As for the
results\textit{:} (\ref{55}) and (\ref{56}), from the \textit{Example 2},
one can say that they are theoretical news and causality related to the
result (\ref{43}). Namely, since $\overset{+\infty }{\underset{k=1}{\sum }}%
k\sin \left( kt\right) =0$ for $t=\frac{\pi }{2}$, that is $\overset{+\infty 
}{\underset{k=0}{\sum }}\left( 2k+1\right) \left( -1\right) ^{k}=0$, it
follows that $\overset{+\infty }{\underset{k=0}{\sum }}2k\left( -1\right)
^{k}=-\overset{+\infty }{\underset{k=0}{\sum }}\left( -1\right) ^{k}=-\frac{1%
}{2}$.

\end{document}